\documentclass[a4paper,draft]{article}
\usepackage{amsmath,amssymb,amsthm,amscd}

\theoremstyle{plain}
\newtheorem{theorem}{Theorem}
\newtheorem{corollary}[theorem]{Corollary}
\theoremstyle{definition}
\newtheorem{remark}{Remark}
\newtheorem{example}{Example}

\title{Classifying Brumer's quintic polynomials \\ by weak Mordell-Weil groups~~\footnotetext{This research was partially supported by the Ministry of Education, Science, Sports and Culture, Grant-in-Aid for Scientific Research (B) No.~19340011, (C) No.~18540048, (C) No.~19540015, and (C) No.~19540057.}}
\author{Masanari \textsc{Kida}, Y\={u}ichi \textsc{Rikuna}, and Atsushi \textsc{Sato}}
\date{}

\newcommand{\bbQ}{\mathbb Q}
\newcommand{\bbZ}{\mathbb Z}
\newcommand{\bbQab}{\bbQ_{a,b}}
\newcommand{\kab}{k_{a,b}}
\newcommand{\bru}[1]{\mathrm{Bru}\!\left({#1}\right)}
\newcommand{\spl}[2]{\mathrm{Spl}_{{#1}}\left({#2}\right)}
\newcommand{\gal}[1]{\mathrm{Gal}\left({#1}\right)}
\newcommand{\cub}[1]{\mathrm{Cub}\!\left({#1}\right)}

\begin{document}

\maketitle

\begin{abstract}
We develop a general classification theory for Brumer's dihedral quintic polynomials by means of Kummer theory arising from certain elliptic curves.
We also give a similar theory for cubic polynomials.
\end{abstract}

\section{Introduction}\label{section:intro}

Brumer's quintic polynomial $\bru{a,b;X}$ is defined by
\[
\bru{a,b;X}:=X^5+(a-3)X^4-(a-b-3)X^3+(a^2-a-2b-1)X^2+bX+a\in\bbQab[X],
\]
where $\bbQab:=\bbQ(a,b)$, with two parameters $a$ and $b$.
It is well-known that the splitting field of the polynomial is of degree $10$ over the base field $\bbQab$ and the Galois group is isomorphic to the dihedral group $\mathcal D_5$ of order $10$.
Moreover, the polynomial is a generic $\mathcal D_5$-polynomial over $\bbQ$ in the sense of \cite{JLY}.
In particular, by specializing the parameters $a$ and $b$, it gives \emph{all} the $\mathcal D_5$-extensions over $\bbQ$.
Although it is known that the number of the parameters of this polynomial is minimal to be generic,
we can observe empirically that many different specializations of parameters give the same splitting field.
Once two specializations are given, it is comparatively easy to determine whether the corresponding Brumer's polynomials have the same splitting field or not (see~\cite[Theorem 3.1]{KRY} for example).
But it is difficult to \emph{find systematically} specializations whose splitting fields are the same as prescribed one's.
A recent progress for this problem was made by Hoshi and Rikuna~\cite{HR}.
They gave a constructive method to transform any given dihedral quintic polynomial into Brumer's form.
Applying this transformation to the Brumer's polynomial $\bru{a,b;X}$ itself, we obtain another Brumer's polynomial with different parameters.
If we apply the transformation recursively, we usually get infinitely many distinct parameters giving the same splitting field.
The explicit procedure is given by the following theorem.
\begin{theorem}[Hoshi and Rikuna]\label{theorem:HR}
If we set $\varphi(a,b):=p(a,b)/q(a,b)$ with
\[
\begin{aligned}
p(a,b) &:= b^4+a(3a+1)(4a-7)b^2-2a(4a^4-4a^3-40a^2+91a-4)b \\
& \qquad +a(a^6+5a^5-81a^4+352a^3-634a^2-65a-1), \\
q(a,b) &:=
4b^3-(a^2-30a+1)b^2-2a(3a+1)(4a-7)b+a(4a^4-4a^3-40a^2+91a-4),
\end{aligned}
\]
then $\bru{a,b;X}$ and $\bru{a,\varphi(a,b);X}$ have the same splitting field over $\bbQab$.
\end{theorem}

Another progress for this problem was due to Kida, Renault, and Yokoyama~\cite{KRY}.
They prove the following theorem.
\begin{theorem}[Kida, Renault, and Yokoyama~\mbox{\cite[Theorem 7.1]{KRY}}]\label{theorem:KRY}
Let $\mathcal E$  be an elliptic curve defined by 
\begin{equation}\label{equation:KRY}
\mathcal E:\  47y^2=4x^3+28x^2+24x+47.
\end{equation}
For any rational point $P=(x(P),y(P))\in\mathcal E(\bbQ)$, if $\bru{1,x(P);X}$ is irreducible over $\bbQ$, then the splitting field of $\bru{1,x(P);X}$ is the same as that of $\bru{1,0;X}$.
Moreover, there are infinitely many $P\in\mathcal E(\bbQ)$ such that $\bru{1,x(P);X}$ is irreducible over $\bbQ$.
\end{theorem}
The Mordell-Weil group $\mathcal E(\bbQ)$ of $\mathcal E$ is a free abelian group of rank $2$ generated by $P_1:=(-1,1)$ and $P_2:=(0,-1)$.
It is shown that, if $P\not\in\langle P_1+3P_2,5P_2\rangle$, then $\bru{1,x(P);X}$ is irreducible.
Thus there are infinitely many pairs of the parameters giving the same splitting field.
This theorem was proven by a rather ad hoc argument using the fact that the splitting field of $\bru{1,0;X}$ is the Hilbert class field over $\bbQ\left(\sqrt{-47}\right)$.

In this paper, we bring together and generalize these results.
As a matter of fact, we develop a general classification theory for Brumer's polynomials by means of Kummer theory of certain elliptic curves.
As a corollary, it appears that there are infinitely many Brumer's polynomials defining the isomorphic field as $\bru{a,b;X}$ for almost every given pair $a,b$.
We also apply this strategy to classify generic cubic polynomials $X^3+bX+a\in\bbQab[X]$.

\section{Elliptic curves associated with Brumer's quintic polynomials}

To state our results, we give some preliminaries.
For the parameters $a$ and $b$, put
\[
d(a,b):=-4b^3+(a^2-30a+1)b^2+2a(3a+1)(4a-7)b-a(4a^4-4a^3-40a^2+91a-4).
\]
It is known that the polynomial discriminant of $\bru{a,b;X}$ with respect to $X$ is $a^2d(a,b)^2$ and that the splitting field $\spl{\bbQab}{\bru{a,b;X}}$ contains a unique quadratic subfield $\kab:=\bbQab(\sqrt{d(a,b)})$.
In the rest of this paper, when we regard $a,b$ as given rational numbers, we drop the subscripts and write $k$ and $\bbQ$ instead of $\kab$ and $\bbQab$, respectively.
If the splitting field of $\bru{\alpha,\beta;X}$ is isomorphic to that of the given $\bru{a,b;X}$, then the quadratic subfields must coinside.
By writing $d:=d(a,b)$, it yields
\begin{equation}\label{equation:d}
du^2=d(\alpha,\beta)
\end{equation}
with some rational number $u\neq 0$.
Further we set
\[
\alpha=a, \qquad (x,y):=(-4d\beta,4d^2u),
\]
and we have an elliptic curve 
\[
\begin{aligned}
E_{a,b} \ : \ y^2 &= x^3+d(a^2-30a+1)x^2 \\ 
& \qquad -8d^2a(3a+1)(4a-7)x-16d^3a(4a^4-4a^3-40a^2+91a-4)
\end{aligned}
\]
defined over $\bbQ$.
This elliptic curve $E_{a,b}$ has a rational point $P_0=(-4db,4d^2)$ coming from the trivial solution $\alpha=a$, $\beta=b$, $u=1$ in \eqref{equation:d}.
Moreover, the order of the point $P_0$ in the Mordell-Weil group $E_{a,b}\left(\bbQab\right)$ is infinite.
This can be observed by specializing, for instance, $a=1$, $b=0$.
This specialization gives $P_0=(0,1)$ on $\mathcal E$, which is isomorphic to $E_{1,0}$, of infinite order.

The curve $E_{a,b}$ has an isogeny $\phi$ of degree $5$ defined over $\bbQ$.
The explicit formula for $\phi$ can be computed by the following \textsf{MAGMA} code (see~\cite{magma}):
\begin{verbatim}
K<a,b>:=FunctionField(Rationals(),2);
c0:=4*a^5-4*a^4-40*a^3+91*a^2-4*a;
c1:=-24*a^3+34*a^2+14*a;
c2:=-a^2+30*a-1;
d:=-4*b^3-c2*b^2-c1*b-c0;
Eab:=EllipticCurve([0,-d*u2,0,4*d^2*u1,-4^2*d^3*u0]);
F:=Factorization(DivisionPolynomial(Eab,5));
T:=SubgroupScheme(Eab,F[1][1]);
Eabs,phi:=IsogenyFromKernel(T);
phis:=DualIsogeny(phi);
\end{verbatim}
The image $E^\ast_{a,b}$ of $\phi$ is also computed as \texttt{Eabs}:
\[
\begin{aligned}
E^\ast_{a,b}\ : \ y^2 &= x^3+(a^2-30a+1)dx^2-8(26a^4-310a^3+327a^2+315a+26)d^2 \\
&~ \qquad +16(68a^6-1120a^5+3804a^4+1760a^3+6929a^2+1380a+68)d^3.
\end{aligned}
\]
Let 
\[
\phi^\ast \ :\ E^\ast_{a,b} \longrightarrow E_{a,b}
\]
be the dual isogeny of $\phi$, which is computed as \texttt{phis}.
The quotient group $E_{a,b}(\bbQ)/\phi^\ast(E^\ast_{a,b}(\bbQ))$ is finite by the weak Mordell-Weil theorem \cite[Chapter VIII \S 1]{Silverman}.
For each rational point $P=(x(P),y(P))\in E_{a,b}(\bbQ)$, we define
\[
\bru{P;X}:=\bru{a,\dfrac{x(P)}{-4d};X}.
\]

We can now state our main theorem.
\begin{theorem}\label{theorem:quintic}
We have the following assertions.
\begin{enumerate}
\item For any rational point $P\in E_{a,b}(\bbQ)$, the Brumer's polynomial $\bru{P;X}$ is reducible over $\bbQ$ if and only if $P\in\phi^\ast(E^\ast_{a,b}(\bbQ))$.
\item There exists a bijection between the following two finite sets
\[
\{\text{subgroup of order }5\text{ in }E_{a,b}(\bbQ)/\phi^\ast(E^\ast_{a,b}(\bbQ))\}
\]
and
\[
\{\spl{\bbQ}{\bru{P;X}}\mid P\in E_{a,b}(\bbQ)\setminus\phi^\ast(E^\ast_{a,b}(\bbQ))\}.
\]
This map is induced by the correspondence $E_{a,b}(\bbQ)\ni P\longmapsto \bru{P;X}$.
\end{enumerate}
\end{theorem}
We will give the proof of this theorem in the next section.
In the course of the proof, it will become apparent that there is a deep structure behind the theorem coming from Kummer theory of the elliptic curve $E_{a,b}$.

As we noted in the above, on the elliptic curve $E_{a,b}$, the point
\begin{equation}\label{equation:P0}
P_0=(-4db,4d^2)
\end{equation}
is of infinite order over $\bbQab$.
Under specialization, $P_0$ remains to have infinite order for most pairs.
Since the point $P_0$ corresponds to the given $\bru{a,b;X}$ and the coset $P_0+\phi^\ast(E^\ast_{a,b}(\bbQ))$ corresponds to its isomorphic class, we obtain the following corollary.
\begin{corollary}
If the Mordell-Weil group $E_{a,b}(\bbQ)$ is infinite, then there are infinitely many $\beta\in\bbQ$ such that $\bru{a,b;X}$ and $\bru{a,\beta;X}$ have the same splitting field over $\bbQ$.
\end{corollary}

\section{Proof of Theorem~\ref{theorem:quintic}}

In this section, we prove Theorem~\ref{theorem:quintic}.
We continue to use the notation in the previous section.

First we want to clarify the roles of two parameters $a$ and $b$.
The first parameter $a$ determines the isomorphism class over $\overline{\bbQ}$ of the elliptic curve $E^\ast_{a,b}$.
In fact, the $j$-invariant of $E^\ast_{a,b}$ is 
\[
\dfrac{(a^4-12a^3+14a^2+12a+1)^3}{a^5(a^2-11a-1)},
\]
which does not depend on $b$.
Let us define a new elliptic curve $E^\ast_a$ defined over $\bbQ$ by
\[
E^\ast_a \ : \ y^2-(a-1)xy-ay=x^3-ax^2.
\]
The curves $E^\ast_{a,b}$ and $E^\ast_a$ are isomorphic over $k=\mathbb Q(\sqrt{d})$ and the isomorphism $f^\ast:\ E^\ast_{a,b}\longrightarrow E^\ast_a$ is given explicitly by
\[
f^\ast \ : \ (x,y)\longmapsto \left(\dfrac{1}{100d}\,x-\dfrac{2(a^2 -5a +1)}{25},\ \dfrac{a-1}{200d}\,x+\dfrac{1}{1000d\sqrt{d}}\,y-\dfrac{2a^3-12a^2-13a-2}{50}\right).
\]
This means that $E^\ast_{a,b}$ is a quadratic twist of $E^\ast_a$ corresponding to the quadratic extension $k/\bbQ$.
The second parameter $b$ determines $d$ and this quadratic twist.

The curve $E^\ast_a$ has a rational point $(0,0)$ of order $5$.
It determines an isogeny $\lambda^\ast:\ E^\ast_a\longrightarrow E_a$.
Here $E_a=E_a^\ast/\langle (0,0)\rangle$ is explicitly given by
\[
E_a  \ : y^2-(a-1)xy-ay=x^3-ax^2-5a(a^2+2a-1)x-a(a^4+10a^3-5a^2+15a-1)
\]
and the isogeny $\lambda^\ast$ maps 
\begin{equation}\label{equation:lambdas}
x\longmapsto\dfrac{x^5-2ax^4+a(a^2+3a-1)x^3-3a^2(a-1)x^2+a^3(a-3)x+a^4}{x^2(x-a)^2}.
\end{equation}
We have an exact sequence induced from $\lambda^\ast$:
\[
\begin{CD}
0 @>>> \ker\lambda^\ast @>>> E^\ast_a @>{\lambda^\ast}>> E_a @>>> 0.
\end{CD}
\]
Taking Galois cohomology as $\gal{\overline{k}/k}$-modules,
we obtain an exact cohomology sequence 
\[
\begin{CD}
E^\ast_a(k) @>{\lambda^\ast}>> E_a(k) @>>> \mathrm{H}^1\left(k, \ker\lambda^\ast(\overline{k})\right).
\end{CD}
\]
Since $\ker\lambda^\ast\left(\overline{k}\right)=\ker\lambda^\ast(\bbQ)$, the Galois group $\gal{\overline{k}/k}$ acts trivially on the kernel.
Therefore it yields an injection into the group of continuous homomorphisms:
\[
E_a(k)/\lambda^\ast(E^\ast_a(k))\lhook\joinrel\longrightarrow\mathrm{Hom}_{\mathrm{cont}}\!\!\left(\gal{\overline{k}/k}, \ker\lambda^\ast(\bbQ)\right).
\]

By the construction of the curves and the morphisms among them, there exists an isomorphism $f:\ E_{a,b}\longrightarrow E_a$ defined over $k$ such that the diagram
\begin{equation}\label{equation:diagram}
\begin{CD}
0 @>>> \ker\phi^\ast @>>> E^\ast_{a,b} @>{\phi^\ast_{/\bbQ}}>> E_{a,b} @>>> 0 \\
@. @. @V{f^\ast_{/k}}VV @VV{f_{/k}}V @. \\
0 @>>> \ker\lambda^\ast @>>> E^\ast_a @>>{\lambda^\ast_{/\bbQ}}> E_a @>>> 0 \\
\end{CD} 
\end{equation}
is commutative with exact raws.
In fact, the isomorphism $f$ can be given explicitly by
\[
f \ : \ (x,y)\longmapsto \left(\dfrac{1}{4d}\,x-2a,\ \frac{a-1}{8d}\,x+\dfrac{1}{8d\sqrt{d}}\,y-\dfrac{a(2a-3)}{2}\right).
\]

Now we shall show that $f$ induces an injection
\begin{equation}\label{equation:injection-ftilde}
\tilde{f} \ : \  E_{a,b}(\bbQ)/\phi^\ast(E^\ast_{a,b}(\bbQ))\lhook\joinrel\longrightarrow E_a(k)/\lambda^\ast(E^\ast_a(k)). 
\end{equation}
Since $E_a$ and $E^\ast_a$ are defined over $\bbQ$, the Galois group $\gal{k/\bbQ}$ acts on $E_a(k)$ and $E^\ast_a(k)$.
Let $E_a(k)\langle 2\rangle$ be the non-$2$-part of $E_a(k)$ which is defined to be a product of the $p$-primary parts ($p\neq 2$).
If we define $\mathsf E_{\pm}$ by 
\[
\mathsf E_+:= \{P\in E_a(k)\langle 2\rangle \mid \sigma P=P\}\qquad\text{and}\qquad \mathsf E_-:= \{P\in E_a (k)\langle 2\rangle \mid \sigma P=-P\},
\]
where $\sigma$ is the generator of $\gal{k/\bbQ}$, then $E_a(k)\langle 2\rangle$ decomposes into the direct sum of these two groups:
\[
E_a(k)\langle 2\rangle=\mathsf E_+\oplus\mathsf E_-.
\]
It is obvious that $\mathsf E_+=E_a(\bbQ)\langle 2\rangle $.
Furthermore, it is easy to show that $\mathsf E_-$ is canonically identified with the non-$2$-part of the group of $\bbQ$-rational points on the quadratic twist of $E_a$ corresponding to the extension $k/\bbQ$ (see, for example, \cite{Kida}).
In our case, it can be identified with $E_{a,b}(\bbQ)\langle 2\rangle$ under $f$.
The same argument can be applied also to $E_a^\ast$ and we obtain
\[
E^\ast_a(k)\langle 2\rangle=\mathsf E^\ast_+\oplus\mathsf E^\ast_-,
\]
where $\mathsf E^\ast_+:=E^\ast_a(\bbQ)\langle 2\rangle$ and $\mathsf E^\ast_-$ can be identified with $E^\ast_{a,b}(\bbQ)\langle 2\rangle$ under $f^\ast$ with the similar definitions of $\mathsf E^\ast_{\pm}$.

Under this preparation, we are now ready to prove \eqref{equation:injection-ftilde}.
Let $R\in E_{a,b}(\bbQ)\langle 2\rangle$ be a rational point satisfying $f(R)\in\lambda^\ast(E^\ast_a(k))\langle 2\rangle$.
We have $f(R)\in\mathsf E_-$.
Since the dual isogeny $\lambda$ of $\lambda^\ast$ is defined over $\bbQ$, it preserves the Galois action.
Hence there is a point $R^\ast\in\mathsf E^-$ such that $f(R)=\lambda^\ast(R^\ast)$.
By the above identification, we can find $P^\ast\in E^\ast_{a,b}(\bbQ)\langle 2\rangle$ such that $f^\ast(P^\ast)=R^\ast$.
From the commutative diagram \eqref{equation:diagram}, it follows that $\phi^\ast(P^\ast)=R$.
This shows that the kernel of the natural map $E_{a,b}(\bbQ)\longrightarrow E_a(k)/\lambda^\ast(E^\ast_a(k))$ agrees with $\phi^\ast(E^\ast_{a,b}(\bbQ))$.
Here we should notice that both of the indices $[E_{a,b}(\bbQ):\phi^\ast(E^\ast_{a,b}(\bbQ))]$ and $[E_a(k):\lambda^\ast(E^\ast_a (k))]$ are powers of $5$.
Therefore the $2$-primary part does not matter at all.
This completes the proof of \eqref{equation:injection-ftilde}.
This also completes the proof of the following theorem.
\begin{theorem}\label{theorem:injection}
There exists an injective homomorphism
\[
E_{a,b}(\bbQ)/\phi^\ast(E^\ast_{a,b}(\bbQ))\lhook\joinrel\longrightarrow\mathrm{Hom}_{\mathrm{cont}}\!\!\left(\gal{\overline{k}/k},\ker\lambda^\ast(\bbQ)\right).
\]
\end{theorem}

By this theorem, every point $P=(x(P),y(P))\in E_{a,b}(\bbQ)$ defines a `Kummer extension' 
\[
K_P:=k\left((\lambda^\ast)^{-1}(f(P))\right)
\]
over $k$.
This field is generated over $k$ by the $x$-coordinate $\eta$ of the point $(\lambda^\ast)^{-1}(f(P))$ because the $y$-coordinate of this point is in $k$.
We now compute the defining polynomial of the extension $K_P/\bbQ$.
Let $N(X)$ be the numerator of $\lambda^\ast(X)$ in \eqref{equation:lambdas} and $D(X)$ the denominator of $\lambda^\ast(X)$ in \eqref{equation:lambdas} respectively.
Then we have
\[
\dfrac{N(\eta)}{D(\eta)}=f(x(P)).
\]
This implies that
\[
B(X):=N(X)-f(x(P))D(X)\in\bbQ[X]
\]
is the defining polynomial of the extension $K_P/\bbQ$.
It is easy to observe that
\begin{equation}\label{equation:lecacheux}
\dfrac{X^5}{a^4}B\left(\dfrac{a}{X}\right)=\bru{P;X}.
\end{equation}
This implies that the field $K_P$ is the splitting field of $\bru{P;X}$ over $\bbQ$.
It is now clear that the group $E_{a,b}(\bbQ)/\phi^\ast(E^\ast_{a,b}(\bbQ))$ in Theorem~\ref{theorem:injection} classifies the isomorphism classes of $\bru{P;X}$ with quadratic subfield $k$.
Therefore the injection in Theorem~\ref{theorem:injection} can be considered as a non-abelian Kummer theory for $\mathcal D_5$-extensions.
This completes the proof of Theorem~\ref{theorem:quintic}.

\begin{remark}
The same observation as~\eqref{equation:lecacheux} was made by Lecacheux~\cite[\S 2]{Lecacheux} from another point of view. 
\end{remark}

\begin{remark}
Hoshi and Rikuna~\cite{HR} remarked that the map $(x,y)\longmapsto\varphi(x,y)$ in Theorem~\ref{theorem:HR} is the multiplication-by-$2$ map of the elliptic curve $d(a,b)y^2=d(a,x)$, which is isomorphic to $E_{a,b}$ over $\bbQab$.
Hence the theorem shows that $\bru{[2^i]P;X}$ ($i=0,1,2,\dots$) have the same splitting field over $\bbQ$, where $[2^i]$ is the multiplication-by-$2^i$ map of $E_{a,b}$.
Note that their construction is based on invariant theory but not on the theory of elliptic curves.
\end{remark}

\section{Examples (quintic cases)}\label{section:example-quintic}
According to our theorems, each subgroup of order $5$ of the finite group $E_{a,b}(\bbQ)/\phi^\ast(E^\ast_{a,b}(\bbQ))$ corresponds to an isomorphic class of some $\bru{a,\beta;X}$.

\begin{example}
We first consider the case where $a=1$ and $b=0$.
This case is relevant to Theorem~\ref{theorem:KRY}.
We have two isogenous curves:
\[
\begin{aligned}
E_{1,0} &: \ y^2=x^3+1316x^2+212064x+78074896, \\
E^\ast_{1,0} &: \ y^2=x^3+1316x^2-6786048x-21410794352.
\end{aligned}
\]
The curve $E_{1,0}$ is obviously isomorphic to the curve $\mathcal E$ defined by~\eqref{equation:KRY} in Theorem~\ref{theorem:KRY}.
Their Mordell-Weil groups are
\[
E_{1,0}(\bbQ)=\langle P_1,P_2\rangle\cong\bbZ^{\oplus 2},\qquad E^\ast_{1,0}(\bbQ)=\langle Q_1,Q_2\rangle\cong\bbZ^{\oplus 2},
\]
where
\[
P_1:=(-188,8836),\quad P_2:=(0,-8836),\quad Q_1:=(10528,1104500),\quad Q_2:=(4653,276125).
\]

The isogeny $\phi^\ast: \ E^\ast_{1,0}\longrightarrow E_{1,0}$ maps as follows:
\[
\phi^\ast(Q_1)=-P_1+2P_2, \quad \phi^\ast(Q_2)=2P_1+P_2.
\]
Therefore the image is 
\[
\phi^\ast(E^\ast_{1,0}(\bbQ))=\langle P_1+3P_2,5P_2\rangle.
\]
Theorem~\ref{theorem:quintic} shows that this subgroup corresponds to reducible Brumer's polynomials.
This group has index $5$ in the full Mordell-Weil group.
Thus there is only one isomorphism class of Brumer's polynomials and the isomorphism class must contain $\bru{1,0; X}$.
From this, Theorem~\ref{theorem:KRY} follows.

For example, we have $\spl{\bbQ}{\bru{1,0;X}}=\spl{\bbQ}{\bru{1,\beta;X}}$ for the following $\beta$'s: 

\[
\begin{array}{c|cccccc}
\text{rational point} & P_1 & P_1-P_2 & P_1+P_2 & 2P_1-P_2 & 2P_1 & 2P_2 \\ \hline
\text{corresponding }\beta & \vbox to 17pt {}-1 & -6 & 41 & \dfrac{47}{25} & -\dfrac{210}{47} & -\dfrac{293}{47} \\ 
\end{array}
\]
\end{example}

\begin{example}
In this example, we consider the case where $a=b=2$.
We have more isomorphism classes in this case.
The isogenous curves are 
\[
\begin{aligned}
E_{2,2} &: \ y^2=x^3+16280x^2-9812992x+41494937600, \\
E^\ast_{2,2} &: \ y^2=x^3+16280x^2+70092800x-30706253824000.
\end{aligned}
\]
The Mordell-Weil groups are free of rank $3$:
\[
E_{2,2}(\bbQ) = \langle P_1,P_2,P_3\rangle\cong \bbZ^{\oplus 3},\qquad E^\ast_{2,2}(\bbQ) = \langle Q_1,Q_2,Q_3\rangle\cong\bbZ^{\oplus 3},
\]
where
\[
\begin{aligned}
P_1 &:= (30784,6658816), & P_2 &:= (-8584,832352 ), & P_3 &:=(-148,-208088),\\
Q_1 &:= (307840,175232000),& Q_2 &:= (34040,5476000),& Q_3 &:=\left(\tfrac{876160}{9},-\tfrac{876160000}{27}\right).
\end{aligned}
\]
The images of the generators under the isogeny $\phi^\ast:\ E^\ast_{2,2}\longrightarrow E_{2,2}$ are
\[
\phi^\ast(Q_1)=-P_1+P_2-P_3, \quad \phi^\ast(Q_2)=2P_1-2P_2-P_3, \quad \phi^\ast(Q_3)=-3P_1-2P_2-P_3.
\]
Thus we obtain 
\[
\phi^\ast(E^\ast_{2,2}(\bbQ))=\langle P_1+4P_2+2P_3,5P_2,5P_3\rangle,
\]
which is of index $25=5^2$ in $E_{2,2}(\bbQ)$.
Hence there are $\frac{5^2-1}{5-1}=6$ subgroups of order $5$ in the quotient group $E_{2,2}(\bbQ)/\phi^\ast(E^\ast_{2,2}(\bbQ))$.
Thus we have the following $6$ isomorphism classes of $\bru{2,\beta;X}$ with $k_{2,2}=\bbQ(\sqrt{-74})$:
\[
\begin{array}{c|cccccc}
\text{coset representative} & P_2 & P_3 & P_2+P_3 & P_2+2P_3 & P_2-2P_3 & P_2-P_3 \\ \hline
\text{corresponding } \beta & \vbox to 17pt {}-\dfrac{29}{4} & -\dfrac{1}{8} & \dfrac{233}{36} & \dfrac{15619}{2500} & \dfrac{40091}{676} & -\dfrac{7}{4} \\
\end{array}
\]
Since $\beta=2$ corresponds to the rational point $(2368, 350464)=-4P_2-P_3+(P_1+4P_2+2 P_3)\equiv P_2-P_3 \pmod{\phi^\ast(E^\ast_{2,2}(\bbQ))}$, we have $\spl{\bbQ}{\bru{2,2;X}}=\spl{\bbQ}{\bru{2,-7/4;X}}$.
Each class contains infinitely many $\beta$ giving the same isomorphism class.
\end{example}

\begin{example}
If $d(a,b)$ is a square of a rational number, then we have $k=\bbQ$ and the splitting field $\spl{\bbQ}{\bru{a,b;X}}$ is a cyclic extension of degree $5$ over $\bbQ$.
This example deals with such a degenerate case.
If this is the case, $E^\ast_{a,b}$ possesses a $\bbQ$-rational point of order $5$ and this rational
point generates the kernel of $\phi^\ast$.
Let $a=1$ and $b=-18$.
Then the isogenous curves are
\[
\begin{aligned}
E_{1,-18} &: \ y^2=x^3-409948x^2+20578452576x-2360098139294192, \\
E^\ast_{1,-18} &: \ y^2=x^3-409948x^2-658510482432x+647219253560911504.
\end{aligned}
\]
The Mordell-Weil group of each curve is a finite group of order $5$:
\[
E_{1,-18}(\bbQ)=\langle P_1\rangle\cong\bbZ/5\bbZ,\qquad E^\ast_{1,-18}(\bbQ)=\langle Q_1\rangle\cong\bbZ/5\bbZ,
\]
where
\[
P_1:=(1054152,857435524)\qquad\text{and}\qquad Q_1:=(-351384,885780500).
\]
By this specialization, the point $P_0\in E_{1,-18}(\bbQ)$ which is defined by~\eqref{equation:P0} falls into the point of finite order.
The point $Q_1$ is sent to the identity point in $E_{1,-18}(\bbQ)$ by $\phi^\ast$.
Thus the quotient is
\[
E_{1,-18}(\bbQ)/\phi^\ast(E^\ast_{1,-18}(\bbQ))\cong\bbZ/5\bbZ.
\]
In particular, for this pair $(a,b)=(1,-18)$, there are finitely many $\beta$ such that $\spl{\bbQ}{\bru{1,-18;X}}=\spl{\bbQ}{\bru{1,\beta;X}}$.
Of course, this does not mean that there are only finitely many $\bru{\alpha,\beta;X}$ having the same splitting field as $\spl{\bbQ}{\bru{1,-18;X}}$.
It is also not true that every degenerate case has this finiteness property.
In fact, if $a=-7$ and $b=-20$, then $d(-7,-20)=5^8$.
We compute $E_{-7,-20}(\bbQ)\cong\bbZ$ and $E^\ast_{-7,-20}(\bbQ)\cong\bbZ\oplus\bbZ/5\bbZ$.
Therefore every coset contains infinitely many points giving the same splitting field.
\end{example}

\section{Classifying cubic polynomials}

In this section, we apply our theory for classifying a generic cubic polynomial
\[
\cub{a,b;X}:=X^3+bX+a\in\bbQab[X].
\]
with two parameters $a$ and $b$.
The polynomial discriminant of $\cub{a,b;X}$ with respect to $X$ is
\[
d(a,b):=-(4b^3+27a^2),
\]
and the extension $\spl{\bbQab}{\cub{a,b;X}}/\bbQab$ has a unique intermediate quadratic field $k_{a,b}:=\bbQab(\sqrt{d})$ with $d:=d(a,b)$.
As before, we drop the subscript and write $k$ and $\bbQ$ instead of $k_{a,b}$ and $\bbQab$, respectively.

For given $a$ and $b$, we associate an elliptic curve
\[
E_{a,b}\ :\ y^2=x^3-432a^2d^3.
\]
It admits a $\bbQ$-isogeny $\phi$ of degree $3$ to the curve
\[
E_{a,b}^\ast \ :\ y^2=x^3+11664a^2d^3.
\]
The explicit expression of $\phi$ is given by
\[
\phi \ :\ (x,y)\longmapsto\left(\dfrac{x^3-1728a^2d^3}{x^2},\dfrac{(x^3+3456a^2d^3)y}{x^3}\right),
\]
and that of the dual isogeny $\phi^\ast:E_{a,b}^\ast\longrightarrow E_{a,b}$ is written by
\[
\phi^\ast \ :\ (x,y)\longmapsto\left(\dfrac{x^3+46656a^2d^3}{9x^2},\dfrac{(-x^3+93312a^2d^3)y}{27x^3}\right).
\]
Let
\[
E_a \ :\ y^2+216ay=x^3-326592a^2\qquad\text{and}\qquad E_a^\ast \ :\ y^2+216ay=x^3
\]
be quadratic twists of $E_{a,b}$ and $E_{a,b}^\ast$ corresponding to the quadratic extension $k/\bbQ$, respectively.
We have isomorphisms defined over $k$:
\[
f \ :\ E_{a,b}\longrightarrow E_a\qquad \text{and}\qquad f^\ast\ :\ E_{a,b}^\ast\longrightarrow E_a^\ast
\]
given explicitly by
\[
f\ :\ (x,y)\longmapsto\left(\dfrac{9}{d}\,x,-\dfrac{27}{d\sqrt{d}}\,y-108a\right)\qquad\text{and}\qquad
f^\ast\ :\ (x,y)\longmapsto\left(\dfrac{1}{d}\,x,\dfrac{1}{d\sqrt{d}}\,y-108a\right).
\]
Let $\lambda^\ast:E_a^\ast\longrightarrow E_a$ be an isogeny of degree $3$ such that the following diagram is commutative with exact rows:
\[
\begin{CD}
0 @>>> \ker\phi^\ast @>>> E^\ast_{a,b} @>{\phi^\ast_{/\bbQ}}>> E_{a,b} @>>> 0 \\
@. @. @V{f^\ast_{/k}}VV @VV{f_{/k}}V @. \\
0 @>>> \ker\lambda^\ast @>>> E^\ast_a @>>{\lambda^\ast_{/\bbQ}}> E_a @>>> 0 \\
\end{CD} 
\]
As in the quintic case, the Kummer sequence attached to this commutative diagram gives our classification theory for the cubic case.
\begin{theorem}\label{theorem:cubic}
There is an injective homomorphism
\[
E_{a,b}(\bbQ)/\phi^\ast(E_{a,b}^\ast(\bbQ))\lhook\joinrel\longrightarrow\mathrm{Hom}_{\mathrm{cont}}\!\!\left(\gal{\overline{k}/k},\ker\lambda^\ast(\bbQ)\right)
\]
sending the class of $P=(x(P),y(P))\in E_{a,b}(\bbQ)$ to the character whose kernel fixes the splitting field of $\cub{a,\frac{x(P)}{-4d};X}$.
\end{theorem}

\section{Examples (cubic cases)}

\begin{example}
We here consider the case where $a=b=1$.
The corresponding cubic polynomial is $\cub{1,1;X}=X^3+X+1$, whose polynomial discriminant is $d=-31$.
The associated elliptic curves are
\[
E_{1,1}\ :\ y^2=x^3+12869712\qquad\text{and}\qquad E^\ast_{1,1}\ :\ y^2=x^3-347482224.
\]
Their Mordell-Weil groups are calculated as follows:
\[
E_{1,1}(\bbQ) = \langle P_1, P_2\rangle\cong\bbZ^{\oplus 2},\qquad E^\ast_{1,1}(\bbQ) = \langle Q_1, Q_2 \rangle\cong\bbZ^{\oplus 2}
\]
where
\[
P_1:=(124,3844),\quad P_2:=(217,-4805),\quad Q_1:=(2232,103788),\quad Q_2:=(3472,-203732).
\]
Then we have
\[
\phi^\ast(Q_1)=-P_2,\qquad \phi^\ast(Q_2)=3P_1,\qquad \phi^\ast(E^\ast_{1,1}(\bbQ))=\langle P_2,3P_1\rangle,
\]
and
\[
E_{1,1}(\bbQ)/\phi^\ast(E^\ast_{1,1}(\bbQ))=\langle P_2+\phi^\ast(E^\ast_{1,1}(\bbQ))\rangle\cong\bbZ/3\bbZ.
\]
This implies that, if $\cub{1,\beta;X}=X^3+\beta X+1\in\bbQ[X]$ is irreducible over $\bbQ$ and $\spl{\bbQ}{X^3+\beta X+1}\supset\bbQ(\sqrt{-31})$, then both $X^3+\beta X+1$ and $X^3+X+1$ have the same splitting field over $\bbQ$.
\end{example}

We can also use Theorem~\ref{theorem:cubic} for a classification of monic cubic polynomials over $\bbQ$ with a fixed polynomial discriminant.

Let $D\neq 0$ be a rational number and $E_D$ an elliptic curve over $\bbQ$ defined by
\[
E_D\ :\ y^2=x^3-432D.
\]
For a rational point $P=(x(P),y(P))\in E_D(\bbQ)$, define
\[
F(P;X):=X^3-\dfrac{x(P)}{12}X-\dfrac{y(P)}{108}\in\bbQ[X].
\]
Then the polynomial discriminant of $F(P;X)$ is $D$.

Let $g(X)\in\bbQ[X]$ be a monic cubic polynomial of discriminant $D$.
Then there exists a $\bbQ$-rational point $P_g=(x(P_g),y(P_g))$ of $E_D$ such that $\spl{\bbQ}{g(X)}=\spl{\bbQ}{F(P_g;X)}$.
Indeed, if we write $g(X)=X^3-pX^2+qX-r$, one can find $P_g=(4(p^2-3q),4(2p^3-9pq+27r))\in E_D(\bbQ)$.
Put
\[
a:=-D^2\qquad\text{and}\qquad b:=-\dfrac{D}{4}\,x(P_g).
\]
Then the elliptic curve $E_D$ is isomorphic to $E_{a,b}$ over $\bbQ$, because $E_{a,b}$ is written as
\[
E_{a,b}\ :\ y^2=x^3-432D\left(\dfrac{D^2}{4}\,y(P_g)\right)^6.
\]
By some computation, we have
\[
\spl{\bbQ}{g(X)}=\spl{\bbQ}{\cub{a,b;X}}=\spl{\bbQ}{\cub{-D^2,-\dfrac{D}{4}\,x(P_g);X}}.
\]
Hence we obtain a classification for monic cubic polynomials of fixed discriminant from Theorem~\ref{theorem:cubic}.

\begin{example}
The image of the weak Mordell-Weil group in $\mathrm{Hom}_{\mathrm{cont}}\!\!\left(\gal{\overline{k}/k}, \ker\lambda^\ast(\bbQ)\right)$ is, in fact, contained in the Selmer group attached to the isogeny $\lambda$.
There are many research on the relationship between the Selmer group and ideal class groups of number fields (see~\cite{Aoki} for the references on this subject).
The following example is one of these phenomena.

Let $D=-3321607$.
This is a fundamental discriminant of $k:=\bbQ(\sqrt{D})$.
It is known that the ideal class group of $k$ is isomorphic to $\bbZ/63\bbZ\oplus (\bbZ/3\bbZ)^{\oplus 2}$.
Hence there are $\frac{3^3-1}{3-1}=13$ sextic fields over $\bbQ$ unramified over $k$ whose Galois groups are isomorphic to $\mathcal D_3$.

The elliptic curve
\[
E_{-3321607}\ : \ y^2=x^3+1434934224
\]
has a $\bbQ$-isogeny $\phi$ of degree $3$ to the curve
\[
E^\ast_{-3321607}\ : \ y^2=x^3-38743224048.
\]
Their Mordell-Weil groups are calculated as follows:
\[
E_{-3321607}(\bbQ)=\langle P_1,\dots, P_6\rangle\cong\bbZ^{\oplus 6}\qquad\text{and}\qquad E^\ast_{-3321607}(\bbQ)=\langle Q_1,\dots, Q_6\rangle\cong\bbZ^{\oplus 6},
\]
where
\[
\begin{aligned}
P_1 &:= (-1124,3860), & P_2 &:= (-1020,19332), & P_3 &:= (-572,35324), \\
P_4 &:= (-167,37819), & P_5 &:= (33,37881), & P_6 &:= (601,40645),\\
Q_1 &:= (3384,2916), & Q_2 &:= (3393,17847), & Q_3 &:= (3897,142965),\\
Q_4 &:= (4608,243108), & Q_5 &:= (5472,353700), & Q_6 &:=(7272,588060).\\
\end{aligned}
\]
The images of the generators under the dual isogeny $\phi^\ast:E_{-3321607}^\ast\longrightarrow E_{-3321607}$ are
\[
\begin{aligned}
\phi^\ast(Q_1) &= 2P_2-P_3-P_4+P_5, & \phi^\ast(Q_2) &= -P_1-2P_2+P_4+2P_5, \\
\phi^\ast(Q_3) &= P_4+P_5-P_6, & \phi^\ast(Q_4) &= P_1-P_2+P_5, \\
\phi^\ast(Q_5) &= P_1-2P_2+2P_3+P_4-2P_5+P_6, & \phi^\ast(Q_6) &= P_1+P_2-P_3-P_6.
\end{aligned}
\]
Hence we have
\[
\phi^\ast(E^\ast_{-3321607}(\bbQ))=\langle P_1+P_3+2P_5+P_6,P_2+P_3+P_5+P_6,P_4+P_5+2P_6,3P_3,3P_5,3P_6\rangle
\]
and
\[
E_{-3321607}(\bbQ)/\phi^\ast(E^\ast_{-3321607}(\bbQ))\cong (\bbZ/3\bbZ)^{\oplus 3}.
\]
Therefore $E_{-3321607}(\bbQ)/\phi^\ast(E^\ast_{-3321607}(\bbQ))$ has the $13$ subgroups of order $3$.
They corresponds to the points
\[
P_3,\quad P_4,\quad P_5,\quad P_3\pm P_4,\quad P_4\pm P_5,\quad P_3\pm P_5,\quad P_3\pm P_4\pm P_5.
\]

By direct calculation, we can check that the Kummer extensions corresponding to the $13$ points above are unramified over $k$.
Thus, if a monic irreducible cubic polynomial $g(X)\in\bbQ[X]$ is of discriminant $-3321607$, the splitting field $\spl{\bbQ}{g(X)}$ is unramified over $\bbQ(\sqrt{-3321607})$.
As we noted, every monic cubic polynomial with discriminant $-3321607$ appears in the thirteen isomorphism classes above.
However, it seems that this assertion does not hold in general.
\end{example}

In conclusion, we make some further remarks.

\begin{remark}
An elliptic curve
\[
\mathcal C_a \ :\ y^2+(a^2+a-1)xy+a^3(a-1)y=x^3+a(a-1)x^2
\]
defined over $\bbQ(a)$ has a $7$-torsion point $(0,0)$.
One can also construct a similar theory for a dihedral septic polynomial starting from this curve.

For example, we have an isogeny $\psi:\mathcal C_a\longrightarrow\mathcal C_a/\langle (0,0)\rangle$ of degree $7$ over $\bbQ(a)$ who maps $x\longmapsto \dfrac{N_a(x)}{D_a(x)}$ with
\[
\begin{aligned}
N_a(x) &:= x^7+2a(a-1)(a+1)x^6-a(a-1)(a^5-7a^4+5a^3-3a^2+2a+1)x^5 \\
& \qquad +a^3(a-1)^2(a^4+13a^3-12a^2+9a-6)x^4+a^4(a-1)^3(a^5+7a^4+8a^3-4a^2-a-1)x^3 \\
& \qquad +a^7(a-1)^4(a+1)(3a^2+5a-3)x^2+a^9(a-1)^5(3a^2+3a-1)x+a^{12}(a-1)^6, \\
D_a(x) &:= x^2(x+a^2(a-1))^2(x+a(a-1))^2.
\end{aligned}
\]
We can show that $N_a(X)-bD_a(X)\in\bbQ(a,b)[X]$ is a dihedral septic polynomial with two parameters $a$ and $b$.
Though it is known that there exists a generic $\mathcal D_7$-polynomial over $\bbQ$ by the affirmative answer for the Noether problem due to Furtw\"{a}ngler~\cite{Fur} (see~\cite[\S 5.5]{JLY}), we do not know if our polynomial is a generic polynomial for $\mathcal D_7$.
Thus we do not know if our scheme classifies all $\mathcal D_7$-extensions over $\bbQ$.
\end{remark}

\begin{remark}
In~\cite{Sato}, the remification properties of the splitting fields of $\bru{a,b;X}$ and $\cub{a,b;X}$ are studied as an anologue of classical Kummer extension.
Their ramification is indeed described in terms of reduction properties of elliptic curves.
This is also one of the advantages of considering these field extensions as `Kummer extensions'.
\end{remark}

\medskip
\begin{flushright}
\begin{tabular}{l}
Masanari \textsc{Kida} \\
Department of Mathematics, \\
University of Electro-Communications \\
Chofu, Tokyo 182--8585, Japan \\
\texttt{kida@sugaku.e-one.uec.ac.jp} \\
\\
Y\={u}ichi \textsc{Rikuna} \\
Department of Applied Mathematics, \\
School of Fundamental Science and Engineering, \\
Waseda University \\
3--4--1 Ohkubo, Shinjuku-ku, Tokyo 169--8555, Japan \\
\texttt{rikuna@moegi.waseda.jp} \\
\\
Atsushi \textsc{Sato} \\
Mathematical Institute, \\
Tohoku University \\
Aoba, Sendai 980--8578, Japan \\
\texttt{atsushi@math.tohoku.ac.jp}
\end{tabular}
\end{flushright}
\end{document}